\documentclass{amsart}

\usepackage{amscd,epsfig, amssymb}
\usepackage[all]{xy}
\xyoption{all}

\newtheorem{theorem}{Theorem}[section]
\newtheorem{lemma}[theorem]{Lemma}
\newtheorem{proposition}[theorem]{Proposition}
\newtheorem{corollary}[theorem]{Corollary}

\theoremstyle{definition}
\newtheorem{definition}[theorem]{Definition}
\newtheorem{example}[theorem]{Example}
\newtheorem{remark}[theorem]{Remark}

\numberwithin{equation}{theorem}

\def\ZZ{{\mathbb Z}}
\def\QQ{{\mathbb Q}}
\def\RR{{\mathbb R}}

\def\PP{{\mathbb P}}
\def\AA{{\mathbb A}}
\def\GG{{\mathbb G}}

\def\shA{{\mathcal A}}
\def\shF{{\mathcal F}}
\def\shM{{\mathcal M}}
\def\shO{{\mathcal O}}
\def\shR{{\mathcal R}}

\def\CDiv{{\rm CDiv}}
\def\Cl{{\rm Cl}}
\def\divis{{\rm div}}
\def\Hom{{\rm Hom}}
\def\id{{\rm id}}
\def\lra{\longrightarrow}
\def\ra{\rightarrow}
\def\Spec{{\rm Spec}}
\def\WDiv{{\rm WDiv}}

\begin{document}
\title[Quotient presentations]{Homogeneous coordinates and quotient 
presentations 
for toric varieties}

% Remove or comment out any unused author tags.
% author one   information
\author[A.~A'Campo-Neuen]{Annette A'Campo-Neuen} 
\address{Fachbereich Mathematik, Johannes Gutenberg-Universit\"at,
  55099 Mainz, Germany}
%\curraddr{}
\email{acampo@enriques.mathematik.uni-mainz.de}
%\thanks{}

% author two information
\author[J.~Hausen]{J\"urgen Hausen} 
\address{Fachbereich Math. und Statistik, Universit\"at Konstanz,
  78457 Konstanz, Germany}
%\curraddr{}
\email{Juergen.Hausen@uni-konstanz.de}
%\thanks{}

% author three information
\author[S.~Schroeer]{Stefan Schr\"oer}
\address{Mathematische Fakult\"at, Ruhr-Universit\"at, 
               44780 Bochum, Germany}
%\curraddr{}
\email{s.schroeer@ruhr-uni-bochum.de}
%\thanks{}

%\date{}

% at present the "communicated by" line appears only in ERA and PROC
%\commby{}

%\dedicatory{}

\begin{abstract}
Generalizing  cones over projective toric varieties, we present
arbitrary toric varieties as  quotients of quasiaffine toric varieties.
Such quotient presentations correspond to groups of Weil divisors
generating the topology. Groups comprising Cartier divisors define  
free quotients, whereas  $\QQ$-Cartier divisors define geometric
quotients.  Each quotient presentation yields homogeneous
coordinates. Using homogeneous coordinates, we express quasicoherent 
sheaves in terms of multigraded modules and describe the set of
morphisms into a toric variety.
\end{abstract}

\maketitle

\section{Introduction}

The projective space 
$\PP^n$ is the quotient of the pointed affine space 
$\AA^{n+1}\setminus 0$ by the diagonal 
$\GG_m$-action. A natural question to ask is whether this generalizes
to other toric varieties. Indeed: Cox \cite{Cox 1995b} and others
showed that each toric variety 
$X$ is the quotient of a smooth quasiaffine toric variety
$\hat{X}$. 

This quasiaffine toric variety  $\hat{X}$ and the corresponding
homogeneous coordinate ring $\Gamma(\hat{X},\shO_{\hat{X}})$, however,
are very large and entail
redundant information. For toric varieties with enough invariant 
  Cartier divisors, Kajiwara
\cite{Kajiwara 1998} found smaller homogeneous coordinate rings.

The goal of this paper is to generalize homogeneous coordinates and to
study them from a geometric  viewpoint.
In our language, homogeneous coordinates correspond 
to quotient presentations. Both the constructions of Cox and 
Kajiwara are quotient
presentations; other examples are cones over
quasiprojective toric varieties.
 Given any particular toric variety, our 
approach provides flexibility in the choice of homogeneous coordinate rings.

Roughly speaking, a \emph{quotient presentation} for a toric variety 
$X$ is a quasiaffine toric variety 
$\hat{X}$, together with an affine surjective toric morphism 
$q\colon \hat{X}\ra X$ such that 
groups of invariant Weil divisors on 
$X$ and 
$\hat{X}$ coincide. The global sections 
$S=\Gamma(\hat{X},\shO_{\hat{X}})$ are the corresponding \emph{homogeneous 
coordinates}
for $X$.

Homogeneous coordinates are useful for various purposes.
For example, Cox \cite{Cox 1995a} described the set of morphism 
$r\colon Y\ra X$ from a scheme $Y$ into a smooth toric variety $X$ in terms of 
homogeneous coordinates. Subsequently, Kajiwara \cite{Kajiwara  1998}
generalized this to toric varieties with
enough effective Cartier divisors.
Using homogeneous coordinates, Brion and Vergne 
\cite{Brion; Vergne 1997} determined  Todd classes  on
simplicial toric varieties. Eisenbud, Mustata and Stillman 
\cite{Eisenbud; Mustata; Stillman 2000} recently applied
homogeneous coordinates to  calculate cohomology
groups of coherent sheaves.

This article is divided into five sections. In the first section, we
define the concept of quotient presentations and give a characterization in
terms of fans. Section~2 contains a description of  quotient
presentations in terms of groups of Weil divisors. Such groups of Weil
divisors are not arbitrary. Rather, they generalize the concept of an
ample invertible sheaf or an ample family of sheaves.

In Section~3, we
relate quotient presentations to geometric invariant theory. Quotient
presentations defined by Cartier or $\QQ$-Cartier divisors are free or
geometric quotients, respectively. Because quotients for group actions
tend to be nonseparated, it is natural  (and requires no extra effort)
to consider nonseparated toric prevarieties as well.

In Section~4, we
shall express quasicoherent sheaves on toric varieties in terms of
multigraded modules over homogeneous coordinate rings. In the last
section, we describe the functor $h_X(Y)=\Hom(Y,X)$ represented by a
toric variety $X$ in terms of sheaf data on $Y$ related to homogeneous
coordinates.

%===========================================================
\section{Quotient presentations}

Throughout we shall work over an arbitrary ground field 
$k$. A \emph{toric variety} is an
equivariant   torus embedding 
$T\subset X$, where $X$ is a separated normal algebraic $k$-variety. 
As usual, $N$ denotes the lattice of 1-parameter subgroups of the torus $T$, 
and $M$ is the dual lattice of characters. 
Recall that toric varieties correspond to finite fans 
$\Delta$ in the lattice $N$. We shall encounter \emph{toric 
prevarieties} as well: These are equivariant  
torus embeddings  as above, but with  
$X$ possibly nonseparated. 

Let $q\colon \hat{X}\ra X$ be a surjective toric morphism of toric
prevarieties. Then we have  a pullback homomorphism
$q^*\colon \CDiv^T(X)\ra\CDiv^{\hat{T}}(\hat{X})$  for invariant
Cartier divisors. There is also a  strict transform
for invariant Weil divisors defined as follows. Let  
$U\subset X$ be the union of all $T$-orbits of codimension $\leq 1$,
and $\hat{U}\subset \hat{X}$ its preimage. Each invariant Weil divisor
on $X$  becomes Cartier on $U$, and the composition 
$$
\WDiv^T(X)=\CDiv^T(U)\stackrel{q^*}{\lra}
\CDiv^{\hat{T}}(\hat{U})\subset\WDiv^{\hat{T}}(\hat{U})\subset
\WDiv^{\hat{T}}(\hat{X})
$$
defines the \emph{strict transform} 
$q^\sharp\colon \WDiv^T(X)\ra \WDiv^{\hat{T}}(\hat{X})$ on the groups
of invariant Weil divisors. Note that $q^{\sharp}$ is injective.

\begin{definition}
\label{quotient presentation}
 A \emph{quotient presentation} for a toric prevariety
$X$ is a quasiaffine toric variety 
$\hat{X}$, together with a surjective affine toric morphism 
$q\colon \hat{X}\ra X$ such that the strict transform 
$q^\sharp\colon \WDiv^T(X)\ra \WDiv^{\hat{T}}(\hat{X})$ is bijective.
\end{definition}

\medskip

This notion is local: Given that $\hat{X}$ is quasiaffine, a toric
morphism $q\colon \hat{X}\ra X$ is a quotient presentation if and only
if for each invariant affine  open subset $U\subset X$ the induced
toric morphism $q^{-1}(U)\ra U$ is a quotient presentation.

\begin{example}
The cones  
$\RR_+(1,0)$, and 
$\RR_+(0,1)$ in the lattice
$\hat{N}=\ZZ^2$ define the  quasiaffine toric variety 
$\hat{X}=\AA^2\setminus 0$.  The projection  
$\ZZ^2\ra \ZZ^2/\ZZ(1,1)$ yields a quotient presentation
$q\colon \AA^2\setminus 0\ra X$ for the projective line 
$X=\PP^1$. We could use the projection onto 
$\ZZ^2/\ZZ(1,-1)$ as well. This defines a quotient presentation 
$q\colon \AA^2\setminus 0\ra X$ for the
affine line 
$X=\AA^1\cup \AA^1$ with origin doubled, which is a nonseparated toric 
prevariety.
\end{example}

\medskip

Here comes a characterization of quotient presentations in terms of fans.
For simplicity, we are content with the separated case.
Suppose that $q\colon \hat{X}\ra X$ is a toric morphism of toric varieties
given by a map of fans $Q\colon (\hat{N},\hat{\Delta})\ra (N,\Delta)$.

\begin{theorem}\label{fan description}
The toric morphism $q\colon \hat{X}\ra X$ is a quotient presentation if and
only if the following conditions hold: 
\begin{enumerate}
\item 
The lattice homomorphism 
$Q\colon \hat{N}\ra N$ has finite cokernel.
\item
The fan 
$\hat{\Delta}$ is a subfan of the fan of faces of a
strongly convex cone 
$\bar{\sigma}\subset \hat{N}_\RR$.
\item
The assignment $\sigma \mapsto Q_{\RR}(\sigma)$ defines bijections
  $\hat{\Delta}^{\max} \to \Delta^{\max}$ and $\hat{\Delta}^{(1)}
  \to \Delta^{(1)}$.
\item
For each primitive lattice vector 
$\hat{v}\in \hat{N}$ generating a ray 
$\hat{\rho}\in \hat{\Delta}$, the image 
$Q(\hat{v})\in N$ is a  primitive lattice vector.
\end{enumerate}
\end{theorem}

 \begin{proof} Suppose the conditions hold. The cone 
  $\bar{\sigma}\subset \hat{N}_\RR$ yields  a  toric open  embedding 
  $\hat{X}\subset X_{\bar{\sigma}}$, hence
  $\hat{X}$ is quasiaffine.

  To see that the map $q \colon  \hat X
  \to X$ is surjective, consider an affine chart
  $X_{\sigma} \subset X$, where $\sigma \in \Delta$ is a maximal cone.
  Since $Q$ induces a bijection of maximal cones,
  there is a $\hat \sigma \in \hat \Delta^{\max}$ such that
  $Q_{\RR}(\hat \sigma) = \sigma$.  Moreover, $Q$ was assumed to have
  a finite cokernel, so $q\colon \hat{T}\to T$ is surjective.
  Since $q$ is equivariant,this implies $X_{\sigma} =
  q(X_{\hat \sigma})$.
 
  To check that the map $q \colon  \hat X \to X$ is affine, 
  keep on considering $X_{\sigma}$. It is easy to see that the inverse
  image of $X_{\sigma}$ is
  \begin{equation}\label{toric inverse image}
   q^{-1}(X_{\sigma}) = \bigcup_{\hat \tau \in \hat \Delta; \;
  Q_{\RR}(\hat \tau) \subset \sigma} X_{\hat \tau}.
  \end{equation}
  Using the bijection $\hat{\Delta}^{(1)} \to \Delta^{(1)}$ we see
  that $Q_{\RR}^{-1}(\sigma)$ contains no element of
  $\hat{\Delta}^{(1)} \setminus \hat \sigma^{(1)}$. Consequently, the
  only cones of $\hat \Delta$ mapped by $Q_{\RR}$ into $\sigma$ are
  the faces of $\hat \sigma$. By the above formula, this means
  $q^{-1}(X_{\sigma}) =  X_{\hat \sigma}$. So we see that $q \colon
  \hat X \to X$ is affine.

  It remains to show that the strict transform is bijective. As to
  this, recall first that
  the invariant prime divisors of $X$ are precisely the closures of
  the $T$-orbits $\Spec k[\rho^\perp\cap M] \subset X$ where $\rho \in
  \Delta^{(1)}$.

  We calculate the strict transform of a
  $T$-stable prime divisor $D \subset X$ corresponding to a ray
  $\rho \in \Delta^{(1)}$. 
  Since $\hat{\Delta}^{(1)} \to \Delta^{(1)}$ is
  bijective, there is a unique ray $\hat \rho \in \hat \Delta^{(1)}$ 
  with $Q_{\RR}(\hat \rho) = \rho$. It follows from 
  \ref{toric inverse image} that $q^{\sharp}(D)$ is a multiple of the
  $\hat{T}$-invariant prime divisor $\hat{D}$ corresponding to
  $\hat{\rho}$. 
  Note that $q^{-1}(X_{\rho})=X_{\hat{\rho}}$. To calculate
  the multiplicity of $\hat{D}$ in $q^{\sharp}(D)$, it suffices to 
  determine the pullback of $D \in\CDiv^T(X_{\rho})$ 
  via $q\colon X_{\hat{\rho}}\to X_{{\rho}}$.

  On the affine chart $X_{\rho}$, every invariant Cartier divisor is
  principal, and if $v$ is the primitive lattice vector in $\rho$
  then the assignment 
  $m \mapsto \langle m,  v \rangle D$ induces
  a natural isomorphism $M/\rho^{\perp} \simeq \CDiv^T(X_{\rho})$.
  Since  we have 
  $$  q^{*}(\divis(\chi^{m})) =
  \divis(\chi^{m} \circ q) = \divis(\chi^{m \circ Q}), $$
  the pullback $q^*\colon \CDiv^T(X_{\rho})\to \CDiv^{\hat{T}}(X_{\hat{\rho}})$
  corresponds to the map 
  $Q^*\colon M/\rho^{\perp}\to \hat{M}/\hat{\rho}^{\perp}$.
  By condition (iv), this map is an isomorphism and hence
  $q^*(D)=\hat{D}$.
 
  Again using bijectivity of $\hat \Delta^{(1)} \to \Delta^{1}$, you conclude
  that the strict transform is bijective.
  Thus the conditions are sufficient. Using similar
  arguments, you see that the conditions are also necessary.
\end{proof}

\begin{example}
Suppose $\sigma\subset\RR^3$ be a strongly convex cone generated by four 
extremal 
rays $\RR_+v_1,\ldots,\RR_+v_4$, defining a fan $\Delta$   in 
$N=\ZZ^3$. Let $\hat{\Delta}$ be the fan of all  faces of the first 
quadrant in $\hat{N}=\ZZ v_1\oplus\ldots\oplus\ZZ v_4$. Then the
canonical surjection $Q\colon \hat{N}\ra N$ gives a quotient
presentation.

\begin{figure}[h]
\begin{center}
\epsfbox{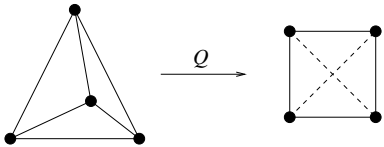}
\caption{A quotient presentation of a non simplicial affine toric
  variety}
\label{nonsimplicial}
\end{center}
\end{figure}

The induced map on nonzero cones looks like Figure~\ref{nonsimplicial}.
The fan $\hat{\Delta}$ comprises 16 cones, whereas  $\Delta$ contains only 10 
cones. You see that two 2-dimensional cones and all 3-dimensional cones in 
$\hat{\Delta}$ map to the maximal cone in $\Delta$. 
\end{example}

\begin{example} 
Let $\Delta$ be a polytopal fan in the lattice $N$, and for each ray
$\rho \in \Delta^{(1)}$ let $v_{\rho} \in \rho$ be the primitive
lattice vector. Consider polytopes $P \subset N_{\RR}$ having edges
$w_{\rho} = n_{\rho}^{-1}v_{\rho}$ with $n_{\rho} \in N$, where
all $\rho\in\Delta^{(1)}$ occur. Each such polytope
defines a quotient presentation of the projective toric variety $X$
associated to $\Delta$:

Set
$\hat{N}=N\oplus\ZZ$. Let $\bar{\sigma}\subset\hat{N}_\RR$ be the cone
generated by $P \times (0,1)$, and $\hat{\Delta}$ the fan of all
strict faces $\hat{\sigma} \subsetneq \bar{\sigma}$. Then the canonical
projection $Q\colon \hat{N}\ra N$ defines a quotient presentation
$q\colon \hat{X}\ra X$. In fact, these quotient presentations are precisely
those obtained from affine cones over $X$. A typical picture is
Figure~\ref{projsurface}.  

\begin{figure}[h]
\begin{center}
\epsfbox{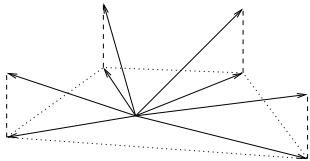}
\caption{A quotient presentation of a projective toric surface}
\label{projsurface}
\end{center}
\end{figure}
\end{example}

%===========================================================
\section{Enough  effective Weil divisors}

The goal of this section is to describe, 
up to isomorphism, the set of all quotient 
presentations of a fixed toric prevariety $X$. 
Recall that we have a canonical map 
$$
\divis\colon M\lra \WDiv^T(X), \quad m\mapsto \divis(\chi^m),
$$
where $\chi^m\in\Gamma(T,\shO_X) $ is the character function corresponding to 
$m\in M $. Suppose $q\colon \hat{X}\ra X$  is a quotient
presentation. The inverse 
$q_* \colon \WDiv^{\hat{T}}(\hat{X})\ra\WDiv^T(X)$ 
of the strict transform yields a factorization 
$$
M\lra \hat{M} \lra \WDiv^T(X)
$$ 
of the canonical map $\divis\colon M\ra\WDiv^T(X)$. We seek to
reconstruct the quotient presentation from such sequences. 

\begin{definition}
A \emph{triangle} is an abstract lattice $\hat{M}$,
together with a sequence 
$M\ra \hat{M} \ra \WDiv^T(X)$, such that the following holds: The composition 
is the   canonical map $\divis\colon M\ra\WDiv^T(X)$, the map
$M\ra \hat{M}$ is injective, and for each invariant affine open subset 
$U\subset X$ 
there is an $\hat{m}\in \hat{M} $ whose image $D\in\WDiv^T(X) $ is effective  
with support $X\setminus U $.
\end{definition}

\medskip

Roughly speaking, the image of $\hat{M} \ra \WDiv^T(X)$ contains enough Weil 
divisors, such that it generates the topology of $X$.
Recall that a scheme $Y$ is separated if the diagonal morphism
$Y\ra Y\times Y$ is a  closed embedding. We say that 
$Y$ is of \emph{affine intersection} if the diagonal is an affine morphism.
In other words, there is an affine open covering $U_i\subset Y$ such that the 
$U_i\cap U_j$ are affine.

\begin{theorem}\label{triangle}
Let $X$ be a toric prevariety of affine intersection.
For each quotient presentation $q\colon \hat{X}\ra X $, the
corresponding sequence   
$M\ra \hat{M} \ra \WDiv^T(X) $ is a triangle. Up to isomorphism, this 
assignment yields a bijection between quotient presentations and triangles.
\end{theorem}

\begin{proof}
Suppose $q\colon \hat{X}\ra X $ is a quotient presentation. Given an 
invariant affine open subset $U\subset X $, the preimage $\hat{U}\subset 
\hat{X}$ is affine as well. There is an effective invariant principal
divisor $\hat{D} \subset \hat{X} $ with support
$\hat{X}\setminus\hat{U}$, because $\hat{X}$ is quasiaffine.
So $D=q_*(\hat{D}) $ is an effective Weil divisor with support 
$X\setminus U$. By construction, $D\in \WDiv^T(X)$ lies in the image
of $\hat{M}$.

Conversely, suppose that $ M\lra \hat{M} \stackrel{\phi}{\lra} \WDiv^T(X)$ 
is a triangle. Set for short $\tilde{M} := \WDiv^T(X)$. Let
$\tilde{N}$ and $\hat{N}$ denote the dual lattices of $\tilde{M}$ and
$\hat{M}$ respectively. Dualizing the triangle, we obtain a sequence 
$$
\tilde{N}  \stackrel{\psi}{\lra}\hat{N} \stackrel{Q}{\lra} N.
$$
For each prime divisor $E \in \tilde{M}$, let $E^{*} \in \tilde{N}$
denote the dual base vector. For every invariant open set $U \subset
X$ we have the submonoid $\tilde{N}_{+}(U)$ generated by the $E^{*}$,
where $E\in\WDiv^T(U)$ is a prime divisor. Let $\hat{\sigma}_{U}
\subset \hat{N}_{\RR}$ be the cone generated by
$\psi(\tilde{N}_{+}(U))$. For example, $\hat{\sigma}_{T} = \{0\}$.

We claim that $\hat{\sigma}_U\subset \hat{\sigma}_X$ is a face
provided that $U\subset X $ is an \emph{affine} invariant open subset.
Indeed by assumption, 
there is an $\hat{m}\in\hat{M} $ such that $D=\phi(\hat{m})$
is an effective Weil divisor with support $X\setminus U$. So for each
prime divisor $E \in \tilde{M}$, we have 
$$
\langle \psi(E^{*}), \hat{m}\rangle = \langle
E^{*},\phi(\hat{m})\rangle = \langle E^{*}, D \rangle \geq 0,
$$
with equality if and only if $E^{*} \in \tilde{N}_+(U)$. So $\hat{m}$
is a supporting hyperplane for $\hat{\sigma}_{X}$ cutting out
$\hat{\sigma}_{U}$ and the claim is verified. In particular, since
$\hat{\sigma}_{T} = \{0\}$ is a face of $\hat{\sigma}_{X}$, this cone
is strictly convex.

For later use, let us also calculate $Q(\psi(E^*))$. 
If $v_{\rho}$ denotes the primitive
lattice vector in the ray $\rho$ corresponding to  
the divisor $E\in\tilde{M}$, we have
$$\langle Q(\psi(E^*)),m\rangle = \langle E^*,\divis\chi^m\rangle=
\langle v_{\rho}, m \rangle\,.$$
That implies $Q(\psi(E^*))=v_{\rho}$.
So in particular, $\psi(E^*)$ is
a primitive lattice vector, and $\psi$ induces a bijection between
the rays of $\tilde{N}_+(X)$ and $\hat\sigma_X$.

Let $\hat{\Delta} $ be the fan in $\hat{N}$ generated by the faces
$\hat{\sigma}_U $, where $U\subset X $ ranges over all invariant affine
open subsets. By construction, this defines a quasiaffine toric
variety $\hat{X}$.  

It remains to construct the quotient presentation $q\colon \hat{X}\ra X
$. First, we do this locally over an invariant affine open subset
$U\subset X $. Let $\sigma_U\subset N_\RR$ be the corresponding cone,
and let $\hat{U}\subset \hat{X} $ be the affine open subset defined by
$\hat{\sigma}_U$.

Clearly, the map $Q\colon\hat{N} \ra N$ has a finite cokernel, since
$M\to \hat{M}$ was is assumed to be injective. We have 
$Q(\hat{\sigma}_U)=\sigma_U$. Moreover, it follows from we saw above
that the map $Q\circ\psi$ induces a bijection between the sets of
primitive lattice vectors generating the rays  of
$\tilde{N}_+(U)$ and $\sigma_U$. Therefore the induced map $v_{\hat{\rho}}\ra
Q(v_{\hat{\rho}}) $ gives a bijection between the primitive lattice
vectors generating the rays in  $\hat{\sigma}_U $ and $\sigma_U $.

By Proposition \ref{fan description}, the associated toric morphism
$\hat{U}\ra U$ is a quotient presentation. To obtain the desired
quotient presentation $q\colon \hat{X}\ra X$, we glue the local
patches. Let $U_1,U_2 \subset X$ be two affine charts. 
The intersection $U := U_1 \cap U_2$ is affine, and the
rays of $\sigma_U$ are in bijection with the invariant prime divisors
in $U_1\cap U_2$. On the other hand, the rays of 
$\hat{\sigma}_1\cap\hat{\sigma}_2$ are the images of the duals to
the prime divisors in $\WDiv^T(U_1\cap U_2)$. This implies
that $Q(\hat{\sigma}_1\cap\hat{\sigma}_2)=\sigma_U$.
\end{proof}

For the following examples,  assume that 
$X$ is a toric variety without nontrivial torus factor.
 Equivalently,  the map 
$M\ra\WDiv^T(X)$ is injective. Such toric varieties are called 
\emph{nondegenerate}.

\begin{example}
Obviously, the factorization 
$M\ra \WDiv^T(X)\stackrel{\id}{\ra}\WDiv^T(X)$ is a triangle. The
corresponding quotient presentation was introduced by Cox \cite{Cox
  1995b}. It is the largest quotient presentation in the sense that
  it dominates all other nondegenerate quotient presentations of $X$.
\end{example}

\begin{example}
Suppose that for each invariant affine open subset $U\subset X $, the 
complement 
$X\setminus U $ is the support of an effective Cartier divisor. Then the 
factorization 
$M\ra \CDiv^T(X)\ra\WDiv^T(X)$ is a triangle.
The corresponding quotient presentation 
$q\colon \hat{X}\ra X$ was studied by Kajiwara \cite{Kajiwara 1998}.  
He says that 
$X$ has \emph{enough   Cartier divisors}.
Note that such toric varieties are
 \emph{divisorial} schemes in the sense of Borelli \cite{Borelli 1963}.
\end{example}

\begin{example}
Suppose 
$X$ is a quasiprojective toric variety. Choose an   ample Cartier divisor 
$D\in\WDiv^T(X) $. Then $M\ra M\oplus \ZZ D\ra \WDiv^T(X) $ is a triangle.
The corresponding quotient presentation 
$q\colon \hat{X}\ra X$ is nothing  but the 
$\GG_m$-bundle obtained from the vector bundle
$L\ra X$ associated to the ample sheaf $\shO_X(D)$.
\end{example}

\medskip
 
Next, we come to existence of quotient presentations:

\begin{proposition}
\label{existence of quotient presentations}
A toric prevariety  admits a quotient presentation if and only if  
it is of affine intersection.
\end{proposition}

\begin{proof}
Suppose $q \colon  \hat{X} \to X$ is a quotient presentation and consider two 
invariant affine charts $X_1$, $X_2$ of $X$. Since $q$ is an affine
toric morphism, the preimages $\hat{X}_{i} := q^{-1}(X_i)$ are invariant
affine charts of $\hat{X}$.

The restriction of $q$ defines a quotient presentation $\hat{X}_1 \cap
\hat{X}_2 \to X_{1} \cap X_{2}$. Since $\hat{X}$ is separated, the
intersection $\hat{X}_1 \cap \hat{X}_2$ is even
affine. Property~\ref{fan description}~ (iii) implies that the image
$X_{1} \cap X_{2} = q(\hat{X}_1 \cap \hat{X}_2)$ is again an affine
toric variety.  

Conversely, let $ X$ be of affine intersection. Choose a splitting 
$M=M'\oplus M''$, where $M'\subset M$ is the kernel of the canonical
map $M\ra \WDiv^T(X)$. It suffices to show that the canonical
factorization 
$$
M \lra M'\oplus\WDiv^T(X)\lra\WDiv^T(X)
$$
is a triangle. Let $U\subset X $ be an invariant affine open
subset. We have to check that the complement $D=X\setminus U $ is a
Weil divisor. For each
invariant affine open subset $V\subset X $, the intersection $U\cap V $ is 
affine, so $V\cap D $ is a Weil divisor. Hence $D $ is a Weil divisor.
\end{proof}

%===========================================================
\section{Free and geometric quotient presentations}
\label{Geometric quotient presentations and principal bundles}

In this section we shall relate quotient presentations to geometric
invariant theory.
Fix a toric prevariety 
$X$, together with a quotient presentation 
$q\colon \hat{X}\ra X$ defined by a triangle
$M\ra \hat{M}\ra \WDiv^T(X)$. Let 
$G\subset \hat{T}$ be the kernel of the induced homomorphism
$\hat{T}\ra T$ of tori. The question is:
In what sense is 
$X$ a quotient of the 
$G$-action on 
$\hat{X}$?

Note that
$G=\Spec k[W]$, such that 
$W=\hat{M}/M$ is the character group of the group scheme 
$G$. Such group schemes are called \emph{diagonalizable}.
The  
$G$-action on 
$\hat{X}$ corresponds to a
$W$-grading on 
$$
 q_*(\shO_{\hat{X}}) = \shR =\bigoplus_{w\in W} \shR_w
$$
for certain   coherent 
$\shO_{X}$-modules 
$\shR_w$. We call them the \emph{weight modules} of the quotient presentation.
To describe the weight modules, consider the commutative diagram 
$$
\xymatrix{
M \ar[r] \ar[d]_{{Q^{*}}} & 
\WDiv^T(X) \ar[r] \ar[d]_{q^{\sharp}}^{\simeq}  & 
\Cl(X) \ar[r] \ar[d]_{q^{\sharp}} &  
0   \\
{\hat{M}} \ar[r] &
{\WDiv}^{\hat{T}}(\hat{X}) \ar[r] 
&{\Cl}(\hat{X}) \ar[r] & 0 .
}
$$
The snake lemma yields a map $W\ra \Cl(X)$. Hence each character 
$w\in W$ gives an isomorphism class of invariant reflexive 
fractional ideals:

\begin{lemma}
\label{weight module}
Each weight module 
$\shR_w$ is an invariant reflexive fractional 
ideal. The isomorphism class $[\shR_w]\in\Cl(X)$ is the image of $-w$.
\end{lemma}

\begin{proof}
First, suppose that the quotient presentation $q\colon \hat{X}\ra X$ is
defined by an inclusion of rings $k[\sigma^\vee\cap M] \subset
k[\hat{\sigma}^\vee\cap \hat{M}]$. The weight module 
$\shR_w\subset\shR$ is  given by the homogeneous component 
$R_w\subset  k[\hat{\sigma}^\vee\cap \hat{M}]$ of degree 
$w\in W$.

Let $v_{\rho}\in N$ and $v_{\hat{\rho}}\in \hat{N}$ be the primitive
lattice vectors generating the rays in $\sigma^{(1)}$ and
$\hat{\sigma}^{(1)}$, respectively. Choose $\hat{m}\in \hat{M}$
representing $w\in W$.
Note that the $\hat T$-invariant Weil divisor
$q_{*}(\divis(\chi^{\hat m}))$ on $\hat X$ is given by the function
$$
\hat{m}\colon \sigma^{(1)}\lra \ZZ,\quad \rho\mapsto \langle \hat{m},
v_{\hat{\rho}} \rangle.
$$
%As explained in \cite{Kempf et al. 1973} p.~27, 
The reflexive fractional ideal $R\subset  k(X)$ over the ring 
$k[\sigma^\vee\cap M]$ corresponding to $-[w]\in\Cl(X)$ is generated by the 
monomials 
$\chi^m\in k[M]$ with
$m\geq  -\hat{m}$   as functions on 
$\sigma^{(1)}$. Obviously, the map 
$\chi^m\mapsto \chi^{Q^*(m)+ \hat{m}}$ induces the desired bijection 
$R\ra R_w$. This is compatible with localization, hence globalizes.
\end{proof}

Suppose a diagonalizable group scheme $G$ acts on a scheme $Y$. An
invariant affine morphism $f\colon Y\ra Z$ with 
$\shO_Z=f_*(\shO_Y)^G$ is called a \emph{good quotient}. Note that this
implies $f(\bigcap W_i)=\bigcap f(W_i)$ for each
family of invariant closed subsets $W_i\subset Y$. Moreover, $f\colon Y\ra
Z$ is a categorial quotient. 

\begin{proposition}
\label{good quotients}
Each quotient presentation $q\colon \hat{X}\ra X$  is a good quotient for the 
$G$-action on 
$\hat{X}$.
\end{proposition}

\begin{proof} 
The problem is local, so we can assume that 
$q\colon \hat{X}\ra X$ is given by an inclusion of rings
$k[\sigma^\vee\cap M] \subset k[\hat{\sigma}^\vee\cap \hat{M}]$. By Lemma 
\ref{weight module}, 
the ring
of invariants 
$k[\hat{\sigma}^\vee\cap \hat{M}]^G$ is nothing but 
$k[\sigma^\vee\cap M]$. 
\end{proof}

Sometimes we can do even better.
Suppose a  diagonalizable  group scheme 
$G$ acts on a scheme 
$Y$. An invariant morphism  $f\colon Y\ra Z$ such that the corresponding morphism
$G\times_Z Y\ra Y\times _Z Y$, $ (g,y)\mapsto (gy,y)  $ is an isomorphism is 
called a \emph{principal homogeneous} $ G$-space. Equivalently, the projection 
$Y\ra Z $ is a principal 
$ G$-bundle in the flat topology (\cite{Milne 1980} III Prop. 4.1).

\begin{proposition}
\label{principal quotient}
The quotient presentation 
$q\colon \hat{X}\ra X$ is a principal homogeneous 
$G$-space if and only if 
$\hat{M}\ra \WDiv^T(X)$ factors through the group of invariant Cartier 
divisors.
\end{proposition}

\begin{proof}
Suppose that 
$\hat{M}$ maps to 
$\CDiv^T(X)$. According to Lemma \ref{weight module}, the homogeneous 
components
in 
$q_*(\shO_{\hat{X}})=\bigoplus_{w\in W}\shR_w$ are invertible.
You easily check that the multiplication maps $\shR_w\otimes\shR_{w'}\ra 
\shR_{w+w'}$ are 
bijective. So by 
\cite{Grothendieck 1970}
Proposition 4.1, the quotient presentation 
$\hat{X}\ra X$ is a principal homogeneous
$G$-space. Hence the condition is sufficient. Reversing the arguments, you see 
that the
condition is necessary as well. 
\end{proof}

\begin{example}
Regular toric prevarieties  
  are factorial, hence their quotient presentations are
  principal homogeneous spaces. Consequently, an arbitrary quotient
  presentation is a principal homogeneous space in codimension 1.
\end{example}

\medskip

For the next result, let us  recall another  
concept from geometric invariant theory.
Suppose a diagonalizable group scheme 
$G$ acts on a scheme 
$Y$. A good quotient $Y\ra Z$ is called a \emph{geometric quotient} 
if it separates the $G$-orbits.

\begin{proposition}
\label{geometric quotient}
Suppose  
$q\colon \hat{X}\ra X$ is a quotient presentation. 
Then 
$X$ is a geometric quotient for the 
$G$-action on 
$\hat{X}$ if and only if 
$\hat{M}\ra\WDiv^T(X)$ factors through the group of invariant 
$\QQ$-Cartier divisors.
\end{proposition}

\begin{proof}
First, we check  sufficiency. Let  
$\hat{M}'\subset \hat{M}$ be  the preimage of the subgroup
$\CDiv^T(X)\subset\WDiv^T(X)$.  The group scheme $H=\Spec
k[\hat{M}/\hat{M}']$ is finite, so its action on 
$\hat{X}$ is automatically closed. Consequently, the quotient 
$\hat{X}'=\hat{X}/H$ is a geometric quotient. You directly see that $\hat{X'}$ 
is quasiaffine. 
Consider the induced toric morphism $q'\colon \hat{X}'\ra X$. The strict transforms 
in
$$
\WDiv^T(X)\lra\WDiv^{\hat{T}'}(\hat{X}')\stackrel{(q')^{\sharp}}{\lra}
\WDiv^{\hat{T}}(\hat{X})
$$ 
are injective, and their composition is bijective. So the map on the right is 
bijective, hence $q'\colon \hat{X}'\ra X$ is another quotient presentation. By 
construction, its triangle $M\ra \hat{M}'\ra \WDiv^T(X)$ factors through 
$\CDiv^T(X)$. According to Proposition \ref{principal quotient}, it is a 
geometric quotient.
So $q\colon \hat{X}\ra X$ is the composition of two geometric quotients, 
hence a 
geometric quotient.

The condition is also necessary. Suppose 
$X$ is a geometric quotient, that means the fibers 
$q^{-1}(x)$ are precisely the 
$G$-orbits. By definition, 
$G$ acts freely on 
$\hat{T}$. By semicontinuity of the fiber dimension, the stabilizers 
$G_{\hat{x}}\subset G$ for $\hat{x}\in \hat{X}$ must be finite. Note that the 
stabilizers are constant along the 
$\hat{T}$-orbits. Hence the stabilizers generate a finite subgroup 
$H\subset G$. 

Set 
$\hat{X}'=\hat{X}/H$. As above, we obtain a  quotient presentation 
$q'\colon \hat{X}'\ra X$. By construction,  
$X$ is a free  geometric quotient for the action of 
$G'=G/H$. Now \cite{Mumford; Fogarty; Kirwan 1993}, Proposition 0.9, 
ensures that   
$q'\colon \hat{X}'\ra X$ is a principal homogeneous 
$G'$-space. By Proposition \ref{principal quotient}, the triangle
$M\ra \hat{M}'\ra \WDiv^T(X)$ factors through $\CDiv^T(X)$. 
This implies
that 
$\hat{M}\ra\WDiv^T(X)$ factors through the group of invariant 
$\QQ$-Cartier divisors.
\end{proof}

\begin{example}
Simplicial toric varieties  
are $\QQ$-factorial, hence their quotient presentations are geometric 
quotients. 
It follows that arbitrary quotient presentations are  geometric quotients in 
codimension 2. 
\end{example}

%===========================================================
\section{Homogeneous coordinates and multigraded modules}

Throughout this section,  fix a toric prevariety $X$ of affine intersection 
 and  choose a quotient
presentation 
$q\colon \hat{X}\ra X$. The goal of this section is to  relate quasicoherent 
$\shO_{X}$-modules to 
multigraded 
modules over homogeneous coordinate rings. This  generalizes the classical 
approach for 
$X=\PP^n$, and results
of   Cox \cite{Cox 1995b}  and   Kajiwara \cite{Kajiwara 1998} as well.

We propose the following definition of homogeneous coordinates.
By assumption, the toric variety
$\hat{X}$ is quasiaffine, so the affine hull
$\bar{X}=\Spec\Gamma(\hat{X},\shO_{\hat{X}})$ is an affine toric variety. 
We call the ring
$S=\Gamma(\hat{X},\shO_{\hat{X}})$ the \emph{homogeneous coordinate ring}
with respect
to the quotient presentation
$q\colon \hat{X}\ra X$.
Let 
$M\subset \hat{M}\ra\WDiv^T(X)$ be its triangle and set 
$W=\hat{M}/M$. 
The action of the diagonalizable group scheme
$G=\Spec k[W]$ on 
$\hat{X}$ induces a 
$G$-action on the affine hull 
$\bar{X}$, which corresponds to a 
$W$-grading 
$S=\bigoplus S_w$.

Suppose 
$F$ is a 
$W$-graded 
$S$-module. Then 
$F$ corresponds to a quasicoherent 
$G$-linearized 
$\shO_{\bar{X}}$-module 
$\shM$. Let 
$i\colon \hat{X}\ra \bar{X}$ be the open inclusion. The restriction 
$i^*(\shM)$ is a 
$G$-linearized quasicoherent 
$\shO_{\hat{X}}$-module. Because 
$q\colon \hat{X}\ra X$ is affine, this corresponds to a 
$W$-grading on
$$
q_*(i^*(\shM)) =\bigoplus_{w\in W} q_*(i^*(\shM))_w.
$$

\begin{definition}
\label{associated sheaf}
The sheaf 
$\tilde{F}=q_*(i^*(\shM))_0$ is called the \emph{associated}  
$\shO_{X}$-module for the  
$W$-graded 
$S$-module 
$F$.
\end{definition}

\medskip

For example, the 
$\shO_{X}$-module associated to 
$S$ is nothing but 
$\tilde{S}=\shO_{X}$. Clearly, 
$F\mapsto \tilde{F}$ is an exact functor from the category of 
$W$-graded 
$S$-modules to the category of quasicoherent 
$\shO_{X}$-modules. 
You easily check that the functor commutes with  direct  limits  and
sends finitely generated 
modules to coherent 
sheaves.

We can pass from quasicoherent sheaves to graded modules as well.
Suppose 
$\shF$ is a quasicoherent 
$\shO_{X}$-module. Decompose 
$q_*(\shO_{X})=\bigoplus_{w\in W} \shR_w$ into weight modules. Then 
$\Gamma_*(\shF)=\bigoplus_{w\in W}\Gamma(\hat{X},\shF\otimes_{\shO_X}\shR_w)$ is a 
$W$-graded 
$S$-module.

\begin{definition}
\label{associated module}
We call 
$\Gamma_*(\shF)$ the 
$W$-graded 
$S$-module \em{associated} to the quasicoherent 
$\shO_{X}$-module 
$\shF$.
\end{definition}

\medskip

For example, 
$\Gamma_*(\shO_{X})=S$. Obviously, 
$\shF\mapsto \Gamma_*(\shF)$ is a functor from the category of quasicoherent 
$\shO_{X}$-modules to the category of 
$W$-graded 
$S$-modules. 

\begin{proposition}
\label{adjunction}
There is a canonical isomorphism 
$\shF\simeq(\Gamma_*(\shF))^\sim$ for each quasicoherent 
$\shO_{X}$-module 
$\shF$.
\end{proposition}

\begin{proof}
By definition, we have 
$(\Gamma_*(\shF))^\sim = (q_*i^*i_*q^*(\shF))_0$. Since 
$i\colon \hat{X}\ra \bar{X}$ is an open embedding, 
$i^*i_*(\shM)\simeq\shM$ holds for each quasicoherent 
$\shO_{\hat{X}}$-module 
$\shM$. This gives 
$(q_*i^*i_*q^*(\shF))_0 \simeq (q_*q^*(\shF))_0$.
Because 
$\shR_0=\shO_{X}$, we have
$(q_*q^*(\shF))_0\simeq \shF$. Consequently,
$\shF\simeq(\Gamma_*(\shF))^\sim$.
\end{proof}

We see that  the functor $F\mapsto \tilde{F}$ 
from graded modules to quasicoherent sheaves is surjective 
on isomorphism classes. It might happen, however, that 
$\tilde{F}=0$ although 
$F\neq 0$.  The next task is   to understand the condition 
$\tilde{F}=0$. To do so, we first have to generalize the classical notions of 
irrelevant ideals and Veronese subrings.

The reduced  closed subset $\bar{X}\setminus \hat{X}$
is  an  invariant closed subset inside the affine toric variety $\bar{X}$.
We call the corresponding $\hat{M}$-homogeneous ideal $S_+\subset S$ the 
\emph{irrelevant ideal}. 
Note that $S_+=S$ holds if and only if $X$ is affine. 

Suppose
$U\subset X$ is an invariant open subset. Let 
$W_U\subset W$ be the subgroup of all 
$w\in W$ such that the corresponding invariant reflexive fractional  ideal
$\shR_w$ is invertible over 
$U\subset X$. Following a standard notation, we call the subring 
$$
S^{(W_U)} = \bigoplus_{w\in W_U} S_w\subset \bigoplus_{w\in W} S_w = S
$$
the \emph{Veronese subring} with respect to 
$W_U\subset W$. Given a 
$W$-graded 
$S$-module $F$, we have the Veronese submodule 
$F^{(W_U)} = \bigoplus_{w\in W_U} F_w$ as well. This is a 
$W_U$-graded 
$S^{(W_U)}$-module.
 
\begin{theorem}
\label{vanishing}
Suppose 
$F$ is a finitely generated
$W$-graded 
$S$-module. Then the condition
$\tilde{F}=0$ holds if and only if   there is an invariant affine open 
covering 
$X=U_1\cup\ldots\cup U_n$ such that some power   
$S_+^k\subset S$ of the irrelevant ideal annihilates the     Veronese 
submodules   
$F^{(W_{U_1})},\ldots, F^{(W_{U_n})}$.
\end{theorem}

\begin{proof} 
Choose 
$\hat{m}_i\in \hat{M}$ such that the homogeneous elements 
$s_i=\chi^{\hat{m}_i}\in S$ define effective Weil divisors with support 
$X\setminus U_i$. Then 
$I=(s_1,\ldots,s_n)\subset S$ has the same radical as the irrelevant ideal 
$S_+\subset S$.

Suppose that 
$\tilde{F}=0$.
Then the restrictions  
$\shF_i=\tilde{F}\mid U_i$ are zero as well. Note that the preimage
$\hat{U}_i= q^{-1}(U_i)$ is affine, with global section ring $S_{s_i}$. 
The   Veronese subring $S_i=( S_{s_i})^{(W_{U_i})}$
defines a factorization 
$$
\hat{U}_i   \lra \Spec( S_i) \lra U_i.  
$$
According to Proposition \ref{weight module}, the map on the right is a 
principal 
bundle for
the action of the diagonalizable group scheme 
$G_i=\Spec k[W_{U_i}]$. 

Setting 
$F_i=(F_{s_i})^{(W_{U_i})} $, we conclude that the 
$S_i$-module
$F_i$ is zero as well. Hence 
$s_i^{k_i}\cdot F_i=0$ for some integer
$k_i>0$, because 
$F_i$ is of finite type. Consequently
$I^k\cdot F_i =0$ with  
$k=\max \left\{ k_i \right\}$.
This shows that the condition is necessary. The converse is similiar.
\end{proof}

%===========================================================
\section{Morphisms into toric varieties}

Throughout this section, fix a toric prevariety 
$X$ of affine intersection. We seek to describe the functor 
$h_X(Y)=\Hom(Y,X)$ represented by 
$X$ in terms of sheaf data on 
$Y$. Here 
$Y$ ranges over the category of 
$k$-schemes.
To do so, choose a quotient presentation 
$q\colon \hat{X}\ra X$. Let 
$S=\Gamma(\hat{X},\shO_{\hat{X}})$ be the homogeneous coordinate ring and set 
$\bar{X}=\Spec(S)$. 

Given a $k$-scheme $Y$, we shall deal with pairs $(\shA,\varphi)$ such
that $\shA$ is a $W$-graded quasicoherent $\shO_Y$-algebra with
$\shA_{0} = \shO_{Y}$, and $\varphi\colon S\otimes\shO_Y\ra\shA$ is a $W$-graded
homomorphism of $\shO_Y$-algebras. For simplicity, we refer to such pairs 
as \emph{$S$-algebras}. An $S$-algebra $(\shA,\varphi)$ yields a diagram 

$$
\xymatrix{%
   {\bar{X} \times Y} \ar[d] %
&  {\Spec(\shA)} \ar[l]_{\Spec(\varphi)} \ar[r]^{p} \ar@{-->}[d] %
& Y \ar@{-->}[d]^{r_{(A,\varphi)}}\\ 
{\bar{X}}           &  {\hat{X}} \ar[l] \ar[r]^{q} & X %
}
$$

The problem is to construct the dashed arrows. For this, we need a
base-point-freeness condition. Recall that the irrelevant ideal
$S_+\subset S$ is the ideal of  the closed subscheme $\bar{X}
\setminus \hat{X}$.

\begin{definition}
\label{base-point-free}
An $S$-algebra $(\shA,\varphi)$ is called \emph{base-point-free} if
for each $y\in Y$ there is an $\hat{M}$-homogeneous $s\in S_+$ such
that the germ $\varphi(s) := \varphi(s \otimes 1) \in\shA_y$ is a
unit.
\end{definition}

\medskip

This is precisely what we need:

\begin{proposition}
\label{morphism}
Each base-point-free $S$-algebra $(\shA,\varphi)$ defines, in a
canonical way, a morphism $r_{(\shA,\varphi)}\colon Y\ra X$.
\end{proposition}

\begin{proof}
First, we claim that $\Spec(\shA)\ra \bar{X}$ factors through the open
subset $\hat{X}\subset\bar{X}$. For $y \in Y$ choose $s \in S_{+}$
such that $\varphi(s)$ is a unit in $\shA_{y}$. Then $\varphi(s)$ is
invertible on a $p$-saturated neighbourhood of $p^{-1}(y) \subset
\Spec(\shA)$. Clearly, this neighbourhood is mapped into $\bar{X}_{s}
\subset \hat{X}$.

According to \cite{Mumford; Fogarty; Kirwan 1993}
Theorem~1.1, the projection $\Spec(\shA)\ra Y$ is a categorical
quotient for the $G$-action defined by the $W$-grading on $\shA$ (here
we use the assumption $\shO_Y=\shA_0$). The composition $\Spec(\shA)\ra
\hat{X}\ra X$ is $G$-invariant. So the universal property of
categorical quotients gives a commutative diagram
\begin{equation}\label{definition of morphism}
\xymatrix{
\Spec(\shA) \ar[r] \ar[d] &  \hat{X} \ar[d] \\
Y \ar[r] & X, }
\end{equation}
which defines the desired morphism $r_{(\shA,\varphi)}\colon Y\ra X$.
\end{proof}

\begin{remark}\label{functorial}
The  assignment $(\shA,\varphi)\mapsto r_{(\shA,\varphi)}$ is
functorial in the following sense: Given a base-point-free $S$-algebra
$(\shA,\varphi)$ on $Y$ and a morphism $f\colon Y'\ra Y$. Then the preimage
$(\shA',\varphi')=(f^*\shA,f^*\varphi)$ is a base-point-free
$S$-algebra on $Y'$, and the corresponding morphisms satisfy
$r_{(\shA',\varphi')}= r_{(\shA,\varphi)}\circ f$.
\end{remark}

\medskip

We call an $\hat{M}$-homogeneous element $s\in S_+$ \emph{saturated},
if $\hat X_{s} = q^{-1}(q(\hat X_{s}))$ holds.
In that case, $X_s := q(\hat X_{s})$ is an affine invariant
open subset with $\Gamma(X_s,\shO_{X}) = S_{(s)}$. Recall that $X$ is
covered by the sets $X_{s}$ with $s \in S_{+}$ saturated. We define
$Y_{\varphi(s)}\subset Y$ to be the (open) subset of all $y\in Y$
where the germ $\varphi(s)\in \shA_y$ is a unit.

\begin{lemma}
\label{preimage}
With the preceding notation, we have $Y_{\varphi(s)} =
r_{(\shA,\varphi)}^{-1}(X_s)$ for each saturated $s \in S_{+}$.
\end{lemma}

\begin{proof}
Let $y \in Y_{\varphi(s)}$. Then $\varphi(s)$ is invertible on a
neighbourhood of the fibre of $\Spec(\shA) \to Y$ over $y$. Looking at
the commutative diagram \ref{definition of morphism}, we see that $s$
is invertible at some point of the fibre of $q \colon  \hat X \to X$
over $x := r_{(\shA,\varphi)}(y)$. Since $s$ is saturated, this means
$x \in X_{s}$. The reverse inclusion is clear by definition.
\end{proof}

Different base-point-free $S$-algebras may define the same
morphism. To overcome this, we need an equivalence relation. Suppose
$(\shA_1,\varphi_1)$ and $(\shA_2,\varphi_2)$ are two base-point-free 
$S$-algebras. Call them \emph{equivalent} if for each saturated $s\in
S_+$, say of degree $w\in W$, the following holds:
\begin{enumerate}
\item The open subsets $Y_{\varphi_i(s)}\subset Y$ coincide for
  $i=1,2$.
\item Over $Y_{\varphi_1(s)}=Y_{\varphi_2(s)}$, the
  $S_{s}^{(w)}$-algebras $\shA_1^{(w)}$ and $\shA_2^{(w)}$ are
  isomorphic.
\end{enumerate}
Here $S_{s}^{(w)} \subset S_{s}$ is the Veronese subring with degrees
in $\ZZ w \subset W$.

\begin{proposition}
\label{equivalent}
Two base-point-free 
$S$-algebras on $Y$ define the same morphism   
$Y\ra X$ if and only if they are equivalent.
\end{proposition}

\begin{proof} 
Suppose that $(\shA_i,\varphi_i)$ are two base-point-free
$S$-algebras, which define two morphisms $r_i\colon Y\ra X$, with
$i=1,2$. First, assume that $r_1=r_2$. Let $s \in S_{+}$ be
saturated. Using Lemma~\ref{preimage}, we infer
$Y_{\varphi_1(s)}=Y_{\varphi_2(s)}$. To check the second condition for
equivalence, note that
$$
\shA_i^{(w)} \vert_{Y_{\varphi_{i}(s)}} = %
\shO_{Y_{\varphi_{i}(s)}}[\varphi_{i}(s),\varphi_{i}(s)^{-1}] %
\quad \hbox{and} \quad %
S_{s}^{(w)} = \Gamma(X_s,\shO_{X})[s,s^{-1}]$$
are Laurent polynomial algebras. So the
map $\varphi_{1}(s)\mapsto \varphi_{2}(s)$ induces the desired
isomorphism.

Conversely, assume that the base-point-free $S$-algebras are
equivalent. Let $s \in S_{+}$ be saturated, and let $w \in W$ be its
degree. Consider the partial quotients
$$\Spec(\shA_{i}) \to \Spec(\shA_{i}^{(w)}) \to Y_{\varphi_{i}(s)} %
\quad \hbox{and} \quad %
\hat{X}_{s} \to \Spec(S_{s}^{(w)}) \to X_{s}$$
Then the isomorphism $\shA_2^{(w)} \to \shA_1^{(w)}$
induces the identity on $Y_{\varphi_{1}(s)} =
Y_{\varphi_{2}(s)}$. Thus the morphism $\Spec(\shA_{1}^{(w)}) \to
\Spec(S_{s}^{(w)})$ induces both, $r_{i} \colon Y_{\varphi_{i}(s)} \to
X_{s}$.
\end{proof}

We come to the main result of this section:

\begin{theorem}
\label{functorial bijection}
The assignment 
$(\shA,\varphi)\mapsto r_{(\shA,\varphi)}$ yields a functorial bijection 
between 
the set of
equivalence classes of base-point-free 
$S$-algebras on 
$Y$ and the set of morphisms
$Y\ra X$.
\end{theorem}

\begin{proof}
In Remark \ref{functorial}, we   already saw that the assignment is
functorial in $Y$. By Proposition~\ref{equivalent}, it is
well-defined on equivalence classes and gives an injection from the
set of equivalence classes to the set of morphisms. It remains to
check that the identity morphism $\id\colon X\ra X$ arises from a
base-point-free $S$-algebra. Indeed: you easily check that 
$\shR=q_*(\shO_{\hat{X}})$, together with the   adjunction map 
$  S\otimes \shO_{X}\ra \shR$ is a base-point-free $S$-algebra defining
the identity on $X$.
\end{proof}

As an application, we generalize the result of Kajiwara in~\cite{Kajiwara 
1998}:

\begin{proposition}
\label{isomorphic}
Suppose the characteristic sequence $M\subset\hat{M}\ra\WDiv^T(X)$
of the quotient presentation $q\colon \hat{X}\ra X$ factors through
the group of Cartier divisors. Then two base-point-free $S$-algebras
define the same morphism into $X$ if and only if they are isomorphic.
\end{proposition}

\begin{proof}
Let $(\shA,\varphi)$ be a base-point-free $S$-algebra on the scheme $Y$ 
defining a 
morphism $r\colon Y\ra X$. Set $\shR= q_*(\shO_{\hat{X}})$. The map
$\Spec(\shA) \ra \hat{X}$ defines a homomorphism 
$\shR \otimes_{\shO_{X}} \shO_Y \ra \shA$. Clearly, it suffices to show
that this map is 
bijective. The problem is local, so we may assume that $X$ is affine, hence
each weight module $\shR_w\subset \shR$ is trivial and  
$S_+=S$ holds. According to Lemma~\ref{preimage}, for each
$\hat{M}$-homogeneous unit $s\in S$, the 
image $\varphi(s)\in\Gamma(Y,\shA)$ is a global unit. Since each weight module 
$\shR_w$ is generated by such a homogeneous unit, we infer that  
$\shR\otimes\shO_Y\ra \shA$
is bijective.
\end{proof}

In general, the homogeneous components of a base-point-free $S$-algebra might 
be noninvertible. However, this does not happen for quotient
presentations that are principal bundles:

\begin{corollary}
Assumptions as in Proposition \ref{isomorphic}. Then each base-point-free 
$S$-algebra 
$(\shA,\varphi)$ has invertible homogeneous components $\shA_w\subset\shA$.
\end{corollary}

\begin{proof}
By assumption, $\shR=q_*(\shO_{\hat{X}})$ has invertible homogeneous
components. 
By the preceding Proposition, each base-point-free $S$-algebra
$(\shA,\varphi)$  
is isomorphic to the preimage
$r_{(\shA,\varphi)}^*(\shR)$.
\end{proof}

%===========================================================


\begin{thebibliography}{ccccc}

\bibitem[1]{Borelli 1963}
{\sc M.~Borelli}:
Divisorial varieties. 
Pac.\ J.\ Math.\ 13, 378--388 (1963). 

\bibitem[2]{Brion; Vergne 1997}
{\sc M.~Brion} and {\sc M.~Vergne}:
An equivariant Riemann--Roch theorem for complete, simplicial toric varieties.
J.\ Reine Angew.\ Math.\ 482, 67--92 (1997).

\bibitem[3]{Cox 1995a}
{\sc D.~Cox}:
The functor of a smooth toric variety.
Tohoku Math.\ J.\ 47, 251--262 (1995).

\bibitem[4]{Cox 1995b}
{\sc D.~Cox}:
The homogeneous coordinate ring of a toric variety.  
J.\ Algebr.\ Geom.\ 4,  17--50 (1995).

\bibitem[5]{Eisenbud; Mustata; Stillman 2000}
{\sc D.~Eisenbud}, {\sc M.~Mustata} and {\sc M.~Stillman}:
Cohomology on toric varieties and local cohomology
with monomial support.\
math.AG/000159.

\bibitem[6]{Grothendieck 1970}
{\sc A.~Grothendieck}:
Groupes diagonisables.
In:
{\sc M.~Demazure} and {\sc A.~Grothendieck} (eds.),
Sch\'emas en groupes II, pp.~1--36.
Lect.\ Notes Math.\ 152.
Springer, Berlin, 1970.

\bibitem[7]{Kajiwara 1998}
{\sc T.~Kajiwara}:
The functor of a toric variety with 
enough invariant effective Cartier divisors.
Tohoku Math.\ J.\  50, 139--157 (1998).

\bibitem[8]{Milne 1980}
{\sc J.~Milne}:
\'Etale cohomology.
Princeton Math. Series 33. 
Princeton Univ. Press, Princeton, 1980.

\bibitem[9]{Mumford; Fogarty; Kirwan 1993}
{\sc D.~Mumford}, {\sc J.~Fogarty} and {\sc F.~Kirwan}:
Geometric invariant theory. Third edition.
Ergeb.\ Math.\  Grenzgebiete (3) 34. Springer, Berlin, 1993.

\bibitem[10]{Seshadri 1972}
{\sc C.~Seshadri}:
Quotient spaces modulo reductive algebraic groups. 
Ann.\ of Math. (2) 95, 511--556  (1972).

\end{thebibliography}
\end{document}